\documentclass[11pt, letter]{amsart} 

\usepackage{amssymb,amscd, epsfig, mathrsfs, xypic, amsmath,amsthm,setspace, color}

\usepackage{ascmac}
 \usepackage{fancybox} 
 
\xyoption{all}

\newtheorem{thm}{Theorem}[section]  
\newtheorem{cor}[thm]{Corollary} 
\newtheorem{lem}[thm]{Lemma}  
\newtheorem{prop}[thm]{Proposition} 
\newtheorem{df-pr}[thm]{Definition-Proposition}

\theoremstyle{definition} 
\newtheorem{defn}[thm]{Definition}
   
\newtheorem{rem}[thm]{Remark}
 
\newtheorem{exm}[thm]{Example}

\newcommand{\wGr}{\operatorname{w\!Gr}}
\newcommand{\aPl}{\operatorname{aPl}}
\newcommand{\aPlt}{\aPl^{\times}}

\newcommand{\wcalS}{\mathcal{S}_w}
\newcommand{\wcalQ}{\mathcal{Q}_w}

\newcommand{\tS}{\widetilde{S}}
\newcommand{\atS}{\widetilde{aS}}
\newcommand{\wtS}{\widetilde{wS}}
\newcommand{\wS}{\mathrm{wS}}

\newcommand{\CC}{{\mathbb C}}

\newcommand{\NN}{{\mathbb N}}

\newcommand{\QQ}{{\mathbb Q}}

\newcommand{\ZZ}{{\mathbb Z}}

\newcommand{\sfd }{{d}}

\newcommand{\sff }{{\mathsf f}}

\newcommand{\sfn }{{n}}

\newcommand{\frakS}{{\mathfrak S}}

\newcommand{\bfx}{{\mathbf x }}

\newcommand{\calE}{{\mathcal E}}
\newcommand{\calF}{{\mathcal F}}

\newcommand{\calP}{{\mathcal P}}
\newcommand{\calQ}{{\mathcal Q}}

\newcommand{\calS}{{\mathcal S}}

\newcommand{\sfT }{{\mathsf T}}

\newcommand{\sfN }{{\mathsf N}}
\newcommand{\sfM }{{\mathsf M}}

\newcommand{\sfK }{{\mathsf K}}

\newcommand{\sfD }{{\mathsf D}}

\newcommand{\surj}{\twoheadrightarrow}

\newcommand{\lan}{{\langle}}
\newcommand{\ran}{{\rangle}}

\newcommand{\inc}{\hookrightarrow}

\newcommand{\Gr}{\operatorname{Gr}}

\newcommand{\wa}{a^w}
\newcommand{\wb}{b^w}

\newcommand{\Lie}{\operatorname{Lie}}

\newcommand{\SL}{{\mbox{SL}}}

\newcommand{\GL}{{\operatorname{GL}}}

\newcommand{\divi}{(1)}

\newcommand{\Mat}{{\operatorname{Mat}}}

\newcommand{\Det}{{\operatorname{Det}}}

\newsavebox{\savepar}

\numberwithin{equation}{section}

\newcounter{labelflag} 
\newcommand{\labelon}{\setcounter{labelflag}{1}}
\newcommand{\Label}[1]{\ifnum\thelabelflag=1\ifmmode
\makebox[0in][l]{\qquad\fbox{\rm#1}} \else
\marginpar{\vspace{0.7\baselineskip} \hspace{-1.1\textwidth}
\fbox{\rm#1}} \fi \fi \label{#1} } \labelon

\title{Schur polynomials and Weighted Grassmannians}
\author{Hiraku Abe}
\address{Osaka City University Advanced Mathematical Institute, 3-3-138 Sugimoto, Sumiyoshi-ku, Osaka 558-8585, Japan}
\email{hirakuabe@globe.ocn.ne.jp} 
\author{Tomoo Matsumura}
\address{Department of Mathematical Sciences, KAIST, 291 Daehak-ro Yuseong-gu Daejeon 305-701, South Korea}
\email{tomoomatsumura@kaist.ac.kr}

\begin{document}

\maketitle

\begin{abstract}
In this paper, we introduce a family of symmetric polynomials by specializing the factorial Schur polynomials. These polynomials represent the weighted Schubert classes of the cohomology of the weighted Grassmannian introduced by Corti-Reid, and we regard these polynomials as analogue of the Schur polynomials. We show that those twisted Schur polynomials are the characters of certain representations. Thus we give an interpretation of the Schubert structure constants of the weighted Grassmannians as the (rational) multiplicities of tensor products of the representations. Furthermore, we derive two types of determinantal formulas for the weighted Schubert classes, in terms of special weighted Schubert classes, and also in terms of Chern classes of tautological orbi-bundles.
\end{abstract}

\section{Introduction}
Let $\calP(\sfd)$ be the set of partitions with at most $\sfd$ rows. For every $\lambda \in \calP(\sfd)$, the Schur function $s_{\lambda}(x)$ is defined as a symmetric polynomial in the variables $(x_1,\cdots, x_{\sfd})$. They form a $\ZZ$-module basis of the algebra  $\ZZ[x]^{\frakS_{\sfd}}$ of symmetric polynomials  in $x$-variables with the coefficients in $\ZZ$. On the other hand, the Grassmannian $\Gr(\sfd,\sfn)$ of complex $\sfd$-planes in $\CC^{\sfn}$ has the distinguished subvarieties, called Schubert varieties,  indexed by the set $\calP(\sfd,\sfn)$ of all partitions contained in the $\sfd \times (\sfn-\sfd)$ rectangle. Their associated cohomology classes $S_{\lambda}, \lambda \in \calP(\sfd,\sfn)$ form a $\ZZ$-module basis of the cohomology $H^*(\Gr(\sfd,\sfn);\ZZ)$. The Schur functions $s_{\lambda}(x)$ represent the Schubert classes $S_{\lambda}$ for the Grassmannian in a sense that there is a surjective ring homomorphism
\begin{align}\label{The surjection}
\ZZ[x]^{\frakS_{\sfd}} \to H^*(\Gr(\sfd,\sfn);\ZZ)
\end{align}
which sends $s_{\lambda}(x)$ to $S_{\lambda}$ if $\lambda \in \calP(\sfd,\sfn)$, or $0$ otherwise. It is worth noting that the above map gives a representation theoretic interpretation to the structure constants with respect to Schubert classes, as we can regard $\ZZ[x]^{\frakS_{\sfd}}$ as the representation ring of the general linear group $\text{GL}_{\sfd}(\CC)$ and the Schur functions $s_{\lambda}(x)$ correspond to the irreducible representations.

The correspondence (\ref{The surjection}) has been generalized in several situations. For example, the equivariant Schubert classes for the Grassmannians are represented by the \emph{factorial Schur functions}, \textit{cf.} \cite{Okounkov, MolevSagan, KnutsonTao}. This equivariant generalization of (\ref{The surjection}) will be the main tool in this paper. Other such examples include the (double/quantum) Schubert polynomials (\cite{FominGelfandPostnikov, LascouxSchutzenberger}) for the (equivariant/quantum) cohomolgoy of full flag varieties and (factorial) Schur $Q$-polynomials (\cite{Ikeda, Ivanov97, Ivanov04}) for (equivariant) cohomology of Lagrangian Grassmannians. One of the advantages of these correspondences is that we can study the structure constants by multiplying actual polynomials.

In this paper, we will introduce and study a twisting of the (factorial) Schur polynomials to generalize the above pictures to the (equivariant) cohomology of the weighted Grassmannians introduced by Corti-Reid \cite{CortiReid}. Below we summarize only non-equivariant results of this paper to avoid complexity, although we build the correspondence for the equivariant cohomology of the weighted Grassmannians first and then derive the non-equivariant one.

Let $w_1,w_2,\cdots$ be an infinite sequence of non-negative integers and $u$ a positive integer. Let $\wGr(\sfd,\sfn)$ be the weighted Grassmannian
introduced in \cite{CortiReid}. Its rational cohomology has a $\QQ$-basis consisting of the weighted Schubert classes $\wS_{\lambda}$, $\lambda \in \calP(\sfd,\sfn)$ and the structure constants of the cohomology ring with respect to this basis are computed
 in \cite{AbeMatsumura}. 
 
For each $\lambda\in \calP(\sfd)$, let $s_{\lambda}(x|a)$ be the factorial Schur function defined by the formula (\ref{fsdef}) below. We consider the polynomial $s^w_{\lambda}(x)$ obtained by specializing $s_{\lambda}(x|a)$ at ${a_i=-(w_i/u)x_{\bar{\emptyset}}}$ for all $i=1,2,\dots,$ 
where $x_{\bar{\emptyset}} = x_1 + \cdots + x_{\sfd}$.
If $w_1=w_2=\cdots=0$, this is nothing but the usual Schur function $s_{\lambda}(x)$. In Proposition \ref{twisted Schur independence}, we show that those polynomials $s^w_{\lambda}(x), \lambda \in \calP(\sfd)$ form a $\QQ$-basis of $\QQ[x]^{\frakS_{\sfd}}$ and represent the weighted Schur classes:

\vspace{0.1in}
\noindent {\bf Theorem A} (Theorem \ref{twisted Schur goes weighted Schubert} below). \ \emph{The map $\QQ[x]^{\frakS_{\sfd}} \to H^*(\wGr(\sfd,\sfn);\QQ)$ defined by sending $s^w_{\lambda}(x)$ to $\wS_{\lambda}$ if $\lambda \in \calP(\sfd,\sfn)$ and $0$ if otherwise, is a surjective ring homomorphism.}
\vspace{0.07in}

Furthermore we prove the following from the definition of $s^w_{\lambda}(x)$.

\vspace{0.1in}
\noindent {\bf Theorem B }(Theorem \ref{twisted Schur and rep} below). \ \emph{Suppose that $u=1$. For any partition $\lambda\in\calP(\sfd)$, we have $s^w_{\lambda}(x)\in\ZZ[x]^{\mathfrak{S}_{\sfd}}$. Moreover, there exists a representation $V^{\lambda}_w$ of $\GL_{\sfd}(\CC)$ such that $s^w_{\lambda}(x)$ is the character $\emph{ch}(V^{\lambda}_w)$ of $V^{\lambda}_w$.}
\vspace{0.1in}

\noindent  This theorem allows us to interpret the weighted Schubert structure constants in terms of representations. Suppose that $u=1$ and the weights are non-increasing, \textit{i.e.} $w_1\geq w_2\geq \cdots$. Then, for each $\lambda,\mu, \nu\in\calP(\sfd)$, there exist non-negative integers $m_{\lambda\mu}\in\ZZ_{\geq1}$ and $m_{\lambda\mu}^{\nu}\in\ZZ_{\geq0}$ such that 
\begin{align*}
 (V^{\lambda}_w \otimes V^{\mu}_w)^{\oplus m_{\lambda\mu}} = \bigoplus_{\nu\in\calP(\sfd)}  (V^{\nu}_w)^{\oplus m_{\lambda\mu}^{\nu}}
 \quad \text{ as representations of $\text{GL}_{\sfd}(\CC)$}.
\end{align*}
By Theorem A and Theorem B we see that 
\[
\wS_{\lambda}\cdot \wS_{\mu} = \sum_{\nu\in\calP(\sfd,\sfn)} \Big(\frac{m_{\lambda\mu}}{m_{\lambda\mu}^{\nu}}\Big) \wS_{\nu}.
\]
Therefore, we can think of the structure constants $\frac{m_{\lambda\mu}}{m_{\lambda\mu}^{\nu}}$ as \emph{rational multiplicities} of $V^{\nu}_w$ in the tensor product $V^{\lambda}_w \otimes V^{\mu}_w$.

To prove Theorem A, we first obtain its equivariant analogue (Proposition \ref{wa}). This also provides an algebraic proof of two deteminantal formulas for the weighted Schubert classes $\wS_{\lambda}$: one is in terms of \textit{special} weighted Schubert classes, and the other is in terms of the Chern classes of the tautological orbifold vector bundles. Let $(k)\in\calP(\sfd)$ be the one row partition with $k$ boxes. Let $\wcalS \inc \calE_w \surj \wcalQ$ be the sequence of the tautological orbifold vector bundles over the weighted Grassmannian defined in Section \ref{Determinantal Formula}. We show

\vspace{0.1in}
\noindent {\bf Theorem C }(Theorem \ref{weighted Giambelli} and \ref{wS in terms specials} below) \emph{For each $\lambda \in \calP(\sfd,\sfn)$, 
\begin{align*}
\wS_{\lambda} 
&= \det \left[\sum_{k=0}^{\lambda_i-i+j} \big(c_1(\calS_w)/u\big)^{k}h_{k}(w_{\lambda_i-i+\sfd+1}, \dots, w_{\sfn}) c_{\lambda_i-i+j-k}(\calQ_w) \right]_{1\leq i\leq j\leq \sfd}\\
&= \det\left[\sum_{r=0}^{j-1} h_r(w_{\lambda_i-i+1+\sfd},\cdots,w_{\lambda_i-i+j-r+\sfd}) \left(\frac{\wS_{\divi}}{-w_{\bar{\emptyset}}}\right)^{r}  \wS_{(\lambda_i-i+j-r)}\right]_{1\leq i,j\leq \sfd}
\end{align*}
where $w_{\bar{\emptyset}}=w_1+\cdots+w_{\sfd}+u$.} 
\vspace{0.1in}

These two formulas coincide in the case of ordinary Grassmannians since the special Schubert classes are the Chern classes of the dual of the tautological bundle of Grassmannians. However this is not the case for the weighted Grassmannians. From Theorem C, it follows that the cohomology of the weighted Grassmannian is generated by the Chern classes of the tautological orbifold bundles. At the end, we give the quotient ring description of the cohomology with generators corresponding to the those Chern classes and their relations.
\subsection*{Acknowledgements}
The authors would like to thank the organizers of ``MSJ Seasonal Institute 2012 Schubert calculus" for providing us an excellent environment for discussions on the topic of this paper. The authors would like to show their gratitude also to Takashi Ikeda and Tatsuya Horiguchi for many useful discussions. The first author is particularly grateful to Takashi Otofuji for many helpful comments. The first author is supported by JSPS Research Fellowships for Young Scientists. The second author is supported by the National Research Foundation of Korea (NRF) grants funded by the Korea government (MEST) (No. 2012-0000795, 2011-0001181). He also would like to express his gratitude to the Algebraic Structure and its Application Research Institute at KAIST for providing him an excellent research environment in 2011-2012.
\section{Preliminary}\label{prelim}
Let $\sfd$ and $\sfn$ be positive integers with $\sfd < \sfn$.  Let $\sfM_{\sfn,\sfd}(\CC)^*$ be the space of $\sfn\times \sfd$ complex matrices of rank $\sfd$. The general linear group $\GL_{\sfd}(\CC)$ and $\GL_{\sfn}(\CC)$ naturally acts on $\sfM_{\sfn,\sfd}(\CC)^*$ by the right and left multiplications respectively. Let 
\[
\aPlt(\sfd, \sfn) := \sfM_{\sfn,\sfd}(\CC)^*/\SL_{\sfd}(\CC)
\]
There is the residual action of $\Det := \GL_{\sfd}(\CC)/\SL_{\sfd}(\CC)$ on $\aPlt(\sfd,\sfn)$ from right since $\SL_{\sfd}(\CC)$ is a normal subgroup of $\GL_{\sfd}(\CC)$. We identify $\Det$ with $\CC^{\times}$ by the determinant map. Let $\sfT := (\CC^{\times})^{\sfn}$ be the diagonal torus in $\GL_{\sfn}(\CC)$ embedded in a standard way. We consider the action of the $(\sfn+1)$-torus $\sfK  := \sfT \times \Det$ on $\aPlt(\sfd,\sfn)$.
\begin{defn}[Corti-Reid \cite{CortiReid}] Let $w_1,\cdots, w_{\sfn}$ be non-negative integers and $u$ a positive integer. Let $\sfD_w:=\CC^{\times}$ and consider the map 
\[
\rho_w:\sfD_{w} \to \sfK ;
\quad t \mapsto (t^{w_\sfn},\cdots, t^{w_{1}}; t^u).
\]
The \emph{weighted Grassmannian}  is the quotient variety $\wGr(\sfd,\sfn):= \aPlt(\sfd,\sfn)/\sfD_w$ with the residual action of $\sfT_{\!w}:= \sfK/\sfD_w$. It is a projective variety with at worst orbifold singularities. The ordinary Grassmannian $\Gr(\sfd,\sfn)$ is obtained by setting $w_1=\cdots = w_{\sfn}=0$ and $u=1$ in the above definition. In this case, we can identify $\sfD_w$ and $\sfT_{\!w}$ with $\Det$ and $\sfT$ respectively. Note that the weights $\{w_i\}$ has the reversed order, compared with the one in \cite{CortiReid} and \cite{AbeMatsumura}. 
\end{defn}
It was shown in \cite{AbeMatsumura} that there are the following ring isomorphisms among the rational equivariant cohomologies, 
\begin{equation}\label{isomorphism}
\xymatrix{
H_{\sfT}^*(\Gr(\sfd,\sfn))  \ar[r]^{\sff^*}& H_{\sfK}^*(\aPlt(\sfd,\sfn)) & \ar[l]_{\sff_w^*} H_{\sfT_{\!w}}^*(\wGr(\sfd,\sfn))
}
\end{equation}
where $\sff^*$ and $\sff_w^*$ are the pullback of the following natural maps between Borel constructions
\[
\xymatrix{
E\sfT\times_{\sfT} \Gr(\sfd,\sfn) &E\sfK\times_{\sfK} \aPlt(\sfd,\sfn) \ar[l]_{\sff} \ar[r]^{\sff_w}& E\sfT_{\!w} \times_{\sfT_{\!w}} \wGr(\sfd,\sfn).
}
\]
The isomorphisms $\sff^*$ and $\sff_w^*$ are actually the isomorphisms of algebras over $H^*(B\sfT)$ and $H^*(B\sfT_{\!w})$ respectively. All the cohomologies treated in this paper are rational coefficient singular cohomology.

Let $\calP(\sfd,\sfn)$ be the set of all partitions $\lambda=(\lambda_1\geq\cdots\geq\lambda_{\sfd})$ that fit inside of the $\sfd \times (\sfn-\sfd)$ rectangle.  For each $\lambda\in\calP(\sfd,\sfn)$, there is a $\sfK$-invariant subvariety $\text{a}\Omega_{\lambda}$ in $\aPlt(\sfd,\sfn)$ which coincides with the pullback of the usual Schubert variety $\Omega_{\lambda}$ in $\Gr(\sfd,\sfn)$ by the quotient map $\aPlt(\sfd,\sfn)\to \Gr(\sfd,\sfn)$. We call the associated equivariant cohomology class $\atS_{\lambda}:=[\text{a}\Omega_{\lambda}]_{\sfK} \in H_{\sfK}^*(\aPlt(\sfd,\sfn))$ the $\sfK$-equivariant Schubert class. This class coincides with the pullback of the usual $\sfT$-equivariant Schubert class $\tS_{\lambda}:=[\Omega_{\lambda}]_{\sfT} \in H_{\sfT}^*(\Gr(\sfd,\sfn))$ along $\sff$. Let $\wtS_{\lambda} \in H_{\sfT_{\!w}}^*(\wGr(\sfd,\sfn))$ be the image of $\atS_{\lambda}$ under the inverse of $\sff_w^*$. It is shown in \cite{AbeMatsumura} that $\wtS_{\lambda}, \lambda \in \calP(\sfd,\sfn)$ form a basis of $H_{\sfT_{\!w}}^*(\wGr(\sfd,\sfn))$ as an $H^*(B\sfT_{\!w})$-module.
\begin{defn}
Let $a_l, l\in \NN$ and $x_i, i= 1,\cdots, \sfd$ be indeterminates. Let $\QQ[a]$ be the ring of polynomials in $a$'s and $\QQ[x]^{\frakS_{\sfd}}$ the ring of symmetric polynomials in $x$'s. Let $\QQ[a][x]^{\frakS_{\sfd}}:=\QQ[a]\otimes_{\QQ} \QQ[x]^{\frakS_{\sfd}}$. Let $\calP(\sfd)$ be the set of partitions with at most $\sfd$ rows. For each $\lambda=(\lambda_1,\dots,\lambda_{\sfd}) \in \calP(\sfd)$, let 
\[
\bar\lambda_i:=\lambda_i + (\sfd + 1 -i), \ \ i=1,\dots, \sfd
\]
 so that the sequence $(\bar\lambda_1,\cdots, \bar\lambda_{\sfd})$ is strictly decreasing. The \emph{factorial Schur polynomial} $s_{\lambda}(x|a)$ is defined as follows \cite{Macdonald92}:
\begin{equation}\label{fsdef}
s_{\lambda}(x|a) = \frac{\det \left[ \prod_{p=1}^{\bar\lambda_i -1}(x_j - a_p)  \right]_{1\leq i,j\leq \sfd}}{\prod_{1 \leq i<j \leq \sfd}  (x_i - x_j)}.
\end{equation}
It is well-known that $\{s_{\lambda}(x|a), \lambda \in \calP(\sfd)\}$ is a basis of $\QQ[a][x]^{\frakS_{\sfd}}$ as a $\QQ[a]$-module.

Let $\{y_1,\dots, y_{\sfn}\}$ be the standard basis of the integral lattice $\Lie (\sfT)_{\ZZ}^*$. Let 
\[
b_i=-y_{\sfn+1-i}, \ \ \ i = 1,\dots, \sfn.
\]
We identify $H^*(B\sfT)$ with $\QQ[b_1,\dots, b_{\sfn}]$ as usual. Consider the projection $\QQ[a] \to \QQ[b_1,\dots, b_{\sfn}]$ by sending $a_i$ to $b_i$ if $i=1,\dots, \sfn$ and $0$ otherwise. This way we regard $H_{\sfT}^*(\Gr(\sfd,\sfn))$ as a $\QQ[a]$-module. Then the following is well-known.
\end{defn}
\begin{thm}[\cite{Okounkov, MolevSagan, KnutsonTao}]\label{surjection to T-coh}
There is a surjective morphism of algebras over $\QQ[a]$
\[
\Phi: \QQ[a][x]^{\frakS_{\sfd}} \to H_{\sfT}^*(\Gr(\sfd,\sfn))
\]
that sends $s_{\lambda}(x|a)$ to $\tS_{\lambda}$ if $\lambda \in \calP(\sfd,\sfn)$, and $0$ if otherwise.
\end{thm}
By composing with the isomorphisms in (\ref{isomorphism}), we obtain
\begin{cor}\label{aPlwGrSurj}
There are surjective ring homomorphisms
\[
\Phi_a: \QQ[a][x]^{\frakS_{\sfd}} \to H^*_{\sfK}(\aPlt(\sfd,\sfn)) \ \ \ \ \ \text{ and } \ \ \ \ \Phi_w: \QQ[a][x]^{\frakS_{\sfd}} \to H^*_{\sfT_{\!w}}(\wGr(\sfd,\sfn))
\]
where $\Phi_a:=f^*\circ \Phi$ and $\Phi_w:=(f_w^*)^{-1}\circ \Phi_a$. The maps $\Phi_a$ and $\Phi_w$ send $s_{\lambda}(x|a)$ to $\atS_{\lambda}$ and $\wtS_{\lambda}$ respectively if $\lambda \in \calP(\sfd,\sfn)$, and $0$ if otherwise. 
\end{cor}
\section{Twisting the module structure}\label{Twisting the module structure}
In this section, we will regard $\QQ[a][x]^{\frakS_{\sfd}}$ as an algebra over a certain polynomial ring in order to regard $\Phi_w$ as a homomorphism of rings over the polynomial ring.

Let $\{z\}$ be the standard basis of $\Lie(\Det)_{\ZZ}^*$ and  $\{\gamma\}$ is the standard basis of $\Lie(\sfD_w)_{\ZZ}^*$. The induced map $\rho_w^*: \Lie(\sfK)^*_{\QQ} \to\Lie(\sfD_w)^*_{\QQ}$ is given by $b_i \mapsto -w_i \gamma$ and $z \mapsto u\gamma$, where $\Lie(\sfK)^*_{\QQ}=\Lie(\sfK)^*_{\ZZ}\otimes_{\ZZ}\QQ$ and $\Lie(\sfD_w)^*_{\QQ} = \Lie(\sfD_w)^*_{\ZZ} \otimes_{\ZZ}\QQ$. We identify the kernel of $\rho_w^*$ with $\Lie(\sfT_{\!w})_{\QQ}^*$. Let  
\begin{equation}
\wb_i:= b_i + (w_i/u) z, \ \ \ i=1,\dots,\sfn,
\end{equation}
then $\{\wb_1,\dots,\wb_{\sfn}\}$ is a basis of $\Lie(\sfT_{\!w})_{\QQ}^*$. Thus we can identify $H^*(B\sfT_{\!w})$ with $\QQ[\wb_1,\dots, \wb_{\sfn}]$.

For each $\lambda \in \calP(\sfd,\sfn)$, let 
\[
\varepsilon^{\lambda}_i = \sfn + 1 - \bar\lambda_i, \ \ \ \ \ \forall i =1,\dots,\sfd. 
\]
Let $e_1,\dots, e_n$ be the standard basis of $\CC^{\sfd}$. Then the $\sfT$-fixed points of $\Gr(\sfd,\sfn)$ are the images $p_{\lambda}$ of the matrices $E_{\lambda}:=[e_{\varepsilon^{\lambda}_1} \cdots e_{\varepsilon^{\lambda}_{\sfd}}] \in \Mat_{\sfn,\sfd}(\CC)^*, \lambda \in \calP(\sfd,\sfn)$. The $\sfT_{\!w}$-fixed points in $\wGr(\sfd,\sfn)$ are also given by the images of the matrices $E_{\lambda}$ in the quotient and denoted by $p^{w}_{\lambda}$. The $\SL_{\sfd}(\CC)$-orbit of $E_{\lambda}$ in $\aPl^{\times}(\sfd,\sfn)$ is denoted by $F_{\lambda}$.

\begin{rem}
In \cite{KnutsonTao, AbeMatsumura}, $\varepsilon^{\lambda}_i$ is denoted by $\lambda_i$ and regarded as the location of $1$'s in $01$ strings. 
For example, for the maximum partition $\lambda=((\sfn-\sfd)^{\sfd}) \in \calP(\sfd,\sfn)$, we have $\bar\lambda=(\sfn,\sfn-1,\dots, \sfn-\sfd+1)$ and $\varepsilon^{\lambda}=(1,2,\dots,\sfd)$. For the empty partition $\emptyset=(0,\dots, 0)$, we have $\bar{\emptyset}=(\sfd,\sfd-1,\dots, 1)$ and $\varepsilon^{\emptyset}=(\sfn+1-\sfd, \dots, \sfn-1,\sfn)$. 
\end{rem}


We extend the diagram (\ref{isomorphism}) with the restrictions to fixed points. 
There is a commutative diagram of pullback maps 
\begin{equation}\label{diagram}
\xymatrix{
H_{\sfT}^*(\Gr(\sfd,\sfn)) \ar[d]_{\sff^*} \ar[r]&H_{\sfT}^*(p_{\lambda})=\QQ[b_1,\cdots,b_{\sfn}]\ar[d]\\
H_{\sfK}^*(\aPlt(\sfd,\sfn)) \ar[r]& H_{\sfK}^*(F_{\lambda}) = \QQ[b_1,\cdots,b_{\sfn},z]/(-b_{\bar\lambda}+z)\\
H_{\sfT_{\!w}}^*(\wGr(\sfd,\sfn)) \ar[r]\ar[u]^{\sff_w^*} & H_{\sfT_{\!w}}^*(p_{\lambda}^w) =\QQ[\wb_1, ..., \wb_{\sfn}]\ar[u]\\
}
\end{equation}
where we denote
\[
b_{\bar\lambda}:=  b_{\bar\lambda_1} + \cdots  + b_{\bar\lambda_{\sfd}}.
\]
The stabilizer $\sfK_{\lambda}$ of points of $F_{\lambda}$ is the kernel of the homomorphism
\[
\sfK \to \CC^{\times}; \ \ \ (t_1,\cdots, t_{\sfn}, s) \mapsto s\cdot t_{\varepsilon^{\lambda}_1} \cdots t_{\varepsilon^{\lambda}_{\sfd}}.
\]
Therefore, since $H_{\sfK}^*(F_{\lambda}) = H^*(B\sfK_{\lambda})$, we can identify  $H_{\sfK}^*(F_{\lambda}) = \QQ[b_1,\cdots,b_{\sfn},z]/(-b_{\bar\lambda}+z)$ as in (\ref{diagram}). The vertical maps on the right are also isomorphisms of rings since the kernel of the composition $\sfK_{\lambda} \to \sfK \to \sfT_{\!w}$ is finite for arbitrary weights $w$'s and $u$.
\begin{lem}
Under the isomorphism $\sff^*: H_{\sfT}^*(\Gr(\sfd,\sfn)) \to H_{\sfK}^*(\aPlt(\sfd,\sfn))$, the preimage of $z=z\cdot 1$ is given by
\[
(\sff^*)^{-1}(z)  = \tS_{\divi} + b_{\bar{\emptyset}}
\]
where $(1)$ is the partition $(1,0,\cdots,0)$. 
\end{lem}
\begin{proof}
If we restrict  $\tS_{\divi} \in H_{\sfT}^*(\Gr(\sfd,\sfn))$ to $p_{\lambda}$, we have  $\tS_{\divi}|_{\lambda} = b_{\bar\lambda} -b_{\bar{\emptyset}}$ by Lemma 3 \cite{KnutsonTao}. Therefore, if we restrict  $\atS_{\divi} \in H_{\sfK}^*(\aPlt(\sfd,\sfn))$ to $F_{\lambda}$, then we have $\atS_{\divi}|_{\lambda} = b_{\bar\lambda} -b_{\bar{\emptyset}}= -b_{\bar{\emptyset}} + z$ by the commutative diagram (\ref{diagram}). Thus $(\sff^*)^{-1}(z)   = \tS_{\divi} + b_{\bar{\emptyset}}$.
\end{proof}
By definition, we have $s_{\divi}(x|a)=x_{\bar{\emptyset}} - a_{\bar{\emptyset}}$. Therefore we have $\Phi_a(x_{\bar{\emptyset}}) = \atS_{\divi} + b_{\bar{\emptyset}}=z$.  The $\QQ[a]$-algebra $\QQ[a][x]^{\frakS_{\sfd}}$ has also a $\QQ[x_{\bar{\emptyset}}, a_1, a_2, \dots]$-algebra structure by the obvious multiplication. This extension of the coefficient ring makes $\Phi_a$ a module homomorphism as follows.
\begin{cor}
The map
\[
\Phi_a: \QQ[a][x]^{\frakS_{\sfd}} \to H_{\sfK}^*(\aPlt(\sfd,\sfn))
\] 
is a homomorphism of algebras over $\QQ[x_{\bar{\emptyset}}, a_1, a_2, \dots]$ with respect to
\[
\QQ[x_{\bar{\emptyset}}, a_1, a_2, \dots] \to \QQ[z, b_1,\dots, b_{\sfn}]; \ \ \ \ x_{\bar{\emptyset}} \mapsto z, \text{ and\ } \  a_l \mapsto \begin{cases} b_l & 1 \leq l \leq \sfn \\ 0 & \sfn < l \end{cases}.
\]
\end{cor}

Now we would like to find a subring of $\QQ[x_{\bar{\emptyset}}, a_1, a_2, \dots]$ that corresponds to $\QQ[\wb_1,\dots, \wb_{\sfn}]$ under $\Phi_w$. Choose an infinite sequence of non-negative integers $w_l, l\in \NN$ where $w_{1}, \dots,  w_{n}$ are the ones chosen to define the weighted Grassmannian. Let
\begin{equation}
\wa_l:= a_l + (w_l/u) x_{\bar{\emptyset}}, \ \ \ \ l\in \NN.
\end{equation}
Then apparently $\{\wa_l, l\in\NN\}$ is a set of algebraically independent variables so that $\QQ[\wa]:= \QQ[\wa_1,\wa_2,\cdots]$ is a polynomial ring. Furthermore it is easy to see that there is a canonical identification as rings $\QQ[a][x]^{\frakS_{\sfd}} \cong \QQ[\wa]\otimes_{\QQ} \QQ[x]^{\frakS_d}=:\QQ[\wa][x]^{\frakS_{\sfd}}$ via $x_i \mapsto x_i$ and $a_l \mapsto \wa_l - (w_l/u)x_{\bar{\emptyset}}$. From the above corollary, we arrive at
\begin{prop}\label{wa}
The surjection 
\[
\Phi_w: \QQ[\wa][x]^{\frakS_{\sfd}} \to H^*_{\sfT_{\!w}}(\wGr(\sfd,\sfn))
\]
is a homomorphism of algebras over $\QQ[\wa]$ with respect to
\[
\QQ[\wa] \to \QQ[\wb_1,\dots, \wb_{\sfn} ]; \ \ \ \  \wa_l \mapsto \begin{cases} \wb_l & 1 \leq l \leq \sfn \\ 0 & \sfn < l. \end{cases}
\]
\end{prop}
\begin{rem}
The set $\{s_{\lambda}(x|a), \lambda \in \calP(\sfd)\}$ of the factorial Schur functions is also a basis of $\QQ[\wa][x]^{\frakS_{\sfd}}$ as a $\QQ[\wa]$-module. Indeed, any non-trivial $\QQ[\wa]$-linear relation among $s_{\lambda}(x|a)$'s will give a non-trivial $\QQ[\wb]$-linear relation among $\wtS_{\lambda}$'s for sufficiently large $\sfn$. Therefore, the claim follows from the $\QQ[\wb]$-linear independency for $\wtS_{\lambda}$'s. It should also be noted that Proposition \ref{wa} and Corollary \ref{aPlwGrSurj} imply that the corresponding structure constants of $\QQ[\wa][x]^{\frakS_{\sfd}}$ (over $\QQ[\wa]$) give the structure constants of $H^*_{\sfT_{\!w}}(\wGr(\sfd,\sfn))$ in terms of the equivariant Schubert classes $\wtS_{\lambda}$'s.
\end{rem}
\begin{rem}\label{tautological T}
Let $[M] \in \Gr(\sfd,\sfn) = M^*_{\sfn,\sfd}(\CC)/\GL_{\sfn}(\CC)$ and denote the subspace of $\CC^{\sfn}$ spanned by the columns of $M$ by $\lan M\ran$. The tautological vector bundle is defined by $\calS = \{ ([M], v) \ |\ [M] \in \Gr(\sfd,\sfn), v \in \lan M \ran\} \subset \Gr(\sfd,\sfn) \times \CC^{\sfn}$ and the action of $\sfT$ is defined by $t \cdot ([M], v) = ([t\cdot M], t\cdot v)$ where $t \in \sfT$ is regarded as a diagonal matrix in $\GL_{\sfn}(\CC)$. The element $\Phi(x_{\bar{\emptyset}})=(\sff^*)^{-1}z$ in $H_{\sfT}^*(\Gr(\sfd,\sfn))$ is $-c^{\sfT}_1(\calS)$ where $c^{\sfT}_1(\calS)$ is  the equivariant first Chern class.
\end{rem}
\section{Vanishing Lemma and Restriction to Fixed Points}
In this section, we obtain the formula for the restriction of  the factorial Schur functions and equivariant Schubert classes to the fixed points in the weighted case.  Let $s^w_{\lambda}(x|\wa)$ be the image of $s_{\lambda}(x|a)$ under the identification $\QQ[a][x]^{\frakS_{\sfd}} \cong \QQ[a^{w}][x]^{\frakS_{\sfd}}$ defined by $a_i=\wa_i - (w_i/u) x_{\bar{\emptyset}}$, \textit{i.e.}
\begin{equation*}
s^w_{\lambda}(x|\wa) :=  s_{\lambda}(x | \wa_1 - (w_1/u) x_{\bar{\emptyset}}, \wa_2 - (w_2/u) x_{\bar{\emptyset}}, \dots ).
\end{equation*}
For each $\mu\in \calP(\sfd)$, let
\begin{equation*}
\wa_{\bar \mu} := \wa_{\bar \mu_1} + \cdots + \wa_{\bar\mu_{\sfd}}, \quad w_{\bar \mu}:=w_{\bar \mu_1} + \cdots + w_{\bar\mu_{\sfd}}+u.
\end{equation*}
Consider the evaluation map $\psi_{\mu}^w : \QQ[\wa][x]^{\frakS_d} \to \QQ[\wa]$ defined by the substition
\begin{equation*}
x_i \mapsto \wa_{\bar\mu_i}  - (w_{\bar\mu_i}/w_{\bar \mu}) \cdot \wa_{\bar\mu} \ \ \text{ for all } i=1,\cdots, \sfd.
\end{equation*}
\begin{lem} \label{lemmasub}
Fix $\mu \in \calP(\sfd)$. Let $c^w_l:=\wa_l - (w_l/w_{\bar \mu})\wa_{\bar\mu}$ for each $l \in \NN$, then 
\[
\psi_{\mu}^w(s_{\lambda}^w(x|\wa)) = s_{\lambda}(c^w_{\bar \mu_1},\dots,c^w_{\bar \mu_{\sfd}}|c^w_1,c^w_2, \dots ).
\]
\end{lem}
\begin{proof}
By a direct computation, we see that $\psi_{\mu}^w(x_{\bar{\emptyset}}) = (u/w_{\bar \mu})\wa_{\bar\mu}$. Therefore
\begin{eqnarray*}
\psi_{\mu}^w(s^w_{\lambda}(x|\wa))
&=&\psi_{\mu}^w(s_{\lambda}(x | \wa_1 - (w_1/u) x_{\bar{\emptyset}}, \wa_2 - (w_2/u) x_{\bar{\emptyset}}, \dots ))\\
&=&s_{\lambda}(c^w_{\bar \mu_1},\dots,c^w_{\bar \mu_{\sfd}}|c^w_1,c^w_2, \dots ).
\end{eqnarray*}
\end{proof}
Let $[\lambda]_{-}:=\{\rho\in\calP(\sfd) \mid \rho\subset\lambda, \ |\overline{\rho} \cap \overline{\lambda}|=\sfd-1\}$.
The following proposition is the generalization of the Vanishing Theorem \cite{MolevSagan, Okounkov}.
\begin{prop} 
For each $\lambda,\mu \in \calP(\sfd)$, we have 
\[
\psi_{\mu}^w(s_{\lambda}^w(x|\wa))=\begin{cases}
0 & \text{ if }  \lambda \not\subset \mu\\
\prod_{\rho \in [\lambda]_{-}}   \left( (w_{\bar \rho}/w_{\bar \lambda})\wa_{\bar\lambda} -\wa_{\bar\rho} \right)   & \text{ if }  \lambda=\mu.
\end{cases}
\]
\end{prop}
\begin{proof}
In \cite{MolevSagan, Okounkov}, it is shown that 
\[
s_{\lambda}(a_{\bar\mu}|a) = \begin{cases}
0 & \text{ if }  \lambda \not\subset \mu \\
\prod_{\rho \in [\lambda]_{-}}   \left(a_{\bar\lambda} -a_{\bar\rho} \right)   & \text{ if }  \lambda=\mu.
\end{cases}
\]
From this and Lemma \ref{lemmasub}, the first part of the case is obvious. The second claim follows from the following computation. Suppose $\lambda=\mu$. We compute
\begin{eqnarray*}
c^w_{\bar\lambda} = c^w_{\bar\lambda_1} + \cdots + c^w_{\bar\lambda_{\sfd}}
=\wa_{\bar\lambda} - \frac{w_{\bar\lambda} - u}{w_{\bar\lambda}} \wa_{\bar\lambda} =  \frac{u}{w_{\bar\lambda}} \wa_{\bar\lambda} \\
\end{eqnarray*}
and
\begin{eqnarray*}
c^w_{\bar\rho} = c^w_{\bar\rho_1} + \cdots + c^w_{\bar\rho_{\sfd}}
= \wa_{\bar\rho} - \frac{w_{\bar\rho}-u}{w_{\bar \lambda}}\wa_{\bar\lambda} = \wa_{\bar\rho} - \frac{w_{\bar\rho}}{w_{\bar \lambda}}\wa_{\bar\lambda} + \frac{u}{w_{\bar \lambda}} \wa_{\bar\lambda} \\
\end{eqnarray*}
Thus
\[
\psi_{\lambda}^w(s_{\lambda}^w(x|\wa)) = s_{\lambda}(c^w_{\bar \lambda_1},\dots,c^w_{\bar \lambda_{\sfd}}|c^w_1,c^w_2, \dots ) = \prod_{\rho \in [\lambda]_{-}}   \left(c^w_{\bar\lambda} -c^w_{\bar\rho} \right) =\prod_{\rho \in [\lambda]_{-}}   \left( \frac{w_{\bar \rho}}{w_{\bar \lambda}}\wa_{\bar\lambda} -\wa_{\bar\rho} \right) .
\]
\end{proof}
We also have the following generalization of Lemma, Section 6, \cite{KnutsonTao}.  
\begin{lem} For each $\lambda, \mu \in \calP(\sfd,\sfn)$,  we have
\[
\wtS_{\lambda}|_{\mu}= \psi_{\mu}^w(s_{\lambda}^w(x|\wa))  \big|_{\wa_l=\wb_l, \forall l\in\NN}.
\]
where $\wb_l:=0$ for all $l>\sfn$.
\end{lem}
\begin{proof}
In (\ref{diagram}), the right vertical isomorphisms send $b_i \in \QQ[b_1,\dots, b_{\sfn}]$ to $\wb_i- (w_i/w_{\bar\mu})\wb_{\bar\mu} \in \QQ[\wb_1,\dots, \wb_{\sfn}]$.  We know from Lemma 6 in \cite{KnutsonTao} that 
\[
\tS_{\lambda}|_{\mu}=s_{\lambda}(  b_{\bar\mu_1}, \cdots, b_{\bar\mu_{\sfd}}  | b_1, \dots, b_{\sfn}, 0, \dots   ) \in \QQ[b_1,\dots,b_{\sfn}].
\]
Therefore, 
\begin{eqnarray*}
\wtS_{\lambda}|_{\mu} &=&\left.\left( \tS_{\lambda}|_{\mu}\right)\right|_{b_i \mapsto \wb_i- (w_i/w_{\bar\mu})\wb_{\bar\mu}}\\
&=& s_{\lambda}(  b_{\bar\mu_1}, \cdots, b_{\bar\mu_{\sfd}}  | b_1, \dots, b_{\sfn}, 0, \dots   )\big|_{b_i \mapsto \wb_i- (w_i/w_{\bar\mu})\wb_{\bar\mu}}\\
&=& \psi_{\mu}^w(s_{\lambda}^w(x|\wa))  \big|_{\wa_l=\wb_l, \forall l\in\NN},
\end{eqnarray*}
where the last equality follows from Lemma  \ref{lemmasub}.
\end{proof}
\section{Twisting Schur polynomials}
Recall that, for given non-negative integers $w_l, l \in \NN$ and a positive integer $u$,  we identify $\QQ[a][x]^{\frakS_{\sfd}} \cong \QQ[\wa][x]^{\frakS_{\sfd}}$ by  sending $ a_l$ to $\wa_l - (w_l/u)x_{\bar{\emptyset}}$. The polynomial $s^w_{\lambda}(x | a^w)$ is nothing but the factorial Schur polynomial $s_{\lambda}(x | a)$ after this parameter change.
\begin{defn} For each $\lambda\in \calP(\sfd)$, define
\[
s^w_{\lambda}(x) := s^w_{\lambda}(x | 0) \ \ \ \in \QQ[x]^{\frakS_{\sfd}}.
\]
Equivalently, we have $s^w_{\lambda}(x) = s_{\lambda}(x | - (w_1/u)x_{\bar{\emptyset}}, - (w_2/u)x_{\bar{\emptyset}}, \dots )$. 
\end{defn}
\begin{exm}\label{example1} Since $s_{(1)}(x|a) = x_{\bar{\emptyset}} - a_{\bar{\emptyset}}$, by the substitution $\wa_l - (w_l/u) x_{\bar{\emptyset}} = a_l$, we have
\[
s^w_{(1)}(x|\wa) = (w_{\bar{\emptyset}}/u)x_{\bar{\emptyset}} - \wa_{\bar{\emptyset}}.
\]
By setting $\wa=0$, we have
\begin{equation}\label{swdiv}
s^w_{(1)}(x) = (w_{\bar{\emptyset}}/u) x_{\bar{\emptyset}}.
\end{equation}
Since $s_{(1,1)}(x|a) = \sum_{1\leq i < j \leq \sfd} (x_i - a_i)(x_j - a_{j-1})$, by the substitution, we have
\[
s_{(1,1)}^w(x|\wa) = \sum_{1\leq i < j \leq \sfd} (x_i - \wa_i + (w_i/u) x_{\bar{\emptyset}})(x_j - \wa_{j-1} + (w_{j-1}/u) x_{\bar{\emptyset}}). 
\]
By setting $\wa=0$, we have
\[
s_{(1,1)}^w(x) = \sum_{1\leq i < j \leq \sfd} (x_i  + (w_i/u) x_{\bar{\emptyset}})(x_j  + (w_{j-1}/u) x_{\bar{\emptyset}}). 
\]
\end{exm}
Since $s^w_{\lambda}(x | a^w), \lambda\in\calP(\sfd)$ form a basis of  $\QQ[\wa][x]^{\frakS_{\sfd}}$ as a $\QQ[\wa]$-module, the polynomials $s^w_{\lambda}(x), \lambda \in \calP(\sfd)$ generate $\QQ[x]^{\frakS_{\sfd}}$ over $\QQ$. For each degree part of $\QQ[x]^{\frakS_{\sfd}}$, the numbers of $s_{\lambda}(x)$ and $s^w_{\lambda}(x)$ of the given degree are the same. So we have the following proposition.
\begin{prop}\label{twisted Schur independence}
The polynomials $s^w_{\lambda}(x), \lambda \in \calP(\sfd)$ form a basis of $\QQ[x]^{\frakS_{\sfd}}$.
\end{prop}
For each $\lambda\in \calP(\sfd,\sfn)$, let $\wS_{\lambda}$ be the corresponding weighted Schubert class in $H^*(\wGr(\sfd,\sfn))$. It is the image of $\wtS_{\lambda}$ under the natural map $H_{\sfT_{\!w}}^*(\wGr(\sfd,\sfn))\to H^*(\wGr(\sfd,\sfn))$. Proposition \ref{wa} implies the following theorem.
\begin{thm}\label{twisted Schur goes weighted Schubert}
The map $\QQ[x]^{\frakS_{\sfd}} \to H^*(\wGr(\sfd,\sfn))$ defined by sending $s^w_{\lambda}(x)$ to $\wS_{\lambda}$ if $\lambda \in \calP(\sfd,\sfn)$ and $0$ if otherwise, is a surjective ring homomorphism.
\end{thm}
We conclude this section with the Pieri rule for the twisted Schur functions $s_{\lambda}^w(x)$. For partitions $\lambda, \lambda'\in\calP(\sfd)$, we denote by $\lambda'\rightarrow\lambda$ the condition that $\lambda'\supset\lambda$ and $|\lambda'/\lambda|=1$ where $\lambda'/\lambda$ is the corresponding skew Young diagram and $|\lambda'/\lambda|$ is the number of the boxes in $\lambda'/\lambda$.
\begin{prop}[The Pieri rule]\label{twisted Pieri}
\begin{eqnarray}
s^w_{\divi}(x) \cdot s^w_{\lambda}(x) &=&  \frac{w_{\bar{\emptyset}}}{w_{{\bar{\lambda}}}}  \sum_{\lambda'\to \lambda}s^w_{\lambda'}(x) \label{wPR}
\end{eqnarray}
\end{prop}
\begin{proof}
It is well-known (\textit{cf.} \cite[p.4434]{MolevSagan}) that
\[
s_{\divi}(x|a) s_{\lambda}(x|a) = (a_{{\bar{\lambda}}} - a_{{\bar{\emptyset}}}) s_{\lambda}(x|a) + \sum_{\lambda' \to \lambda} s_{\lambda'}(x|a).
\]
By substituting $a_i \mapsto -(w_i/u)x_{\bar{\emptyset}}$, we have
\[
a_{{\bar{\lambda}}} - a_{{\bar{\emptyset}}} \mapsto 
(x_{\bar{\emptyset}}/u)(- w_{ {\bar{\lambda}}} + w_{ {\bar{\emptyset}}}).
\]
Thus we obtain
\begin{align*}
 \sum_{\lambda' \to \lambda} s^w_{\lambda'}(x) 
 &= s^w_{\divi}(x) s^w_{\lambda}(x) - \frac{- w_{ {\bar{\lambda}}} + w_{ {\bar{\emptyset}}}}{u} x_{\bar{\emptyset}} s^w_{\lambda}(x) \\
 &= s^w_{\divi}(x) s^w_{\lambda}(x)- \frac{- w_{ {\bar{\lambda}}} + w_{ {\bar{\emptyset}}}}{w_{{\bar{\emptyset}}}} s^w_{\divi}(x) s^w_{\lambda}(x) =  \frac{w_{{\bar{\lambda}}}}{w_{\bar{\emptyset}}} s^w_{\divi}(x) s^w_{\lambda}(x).
\end{align*}
where the second equality follows from (\ref{swdiv}). This proves the claim.
\end{proof}
\section{Representations associated to $s^w_{\lambda}(x)$}
In this section, we will see that the twisted Schur function $s^w_{\lambda}(x)$ is  a character of a representation of $\text{GL}_{\sfd}(\CC)$ under the assumption that $u=1$. It basically follows from \cite{MolevSagan}. 

Recall that the character $\text{ch}(V^{\lambda})$ of the irreducible $\text{GL}_{\sfd}(\CC)$-representation $V^{\lambda}$ with the highest weight $\lambda=(\lambda_1,\dots,\lambda_{\sfd})$ coincides with the Schur function $s_{\lambda}(x)$. This gives us the ring isomorphism
\begin{align}\label{rep ring and symm poly}
\text{ch}:  R(\text{GL}_{\sfd}(\CC)) \stackrel{\cong}{\longrightarrow} \ZZ[x]^{\mathfrak{S}_{\sfd}}
\end{align}
where $R(\text{GL}_{\sfd}(\CC))$ is the polynomial representation ring. 
\begin{thm}\label{twisted Schur and rep}
Suppose that $u=1$. For each $\lambda\in\calP(\sfd)$, we have $s^w_{\lambda}(x)\in\ZZ[x]^{\mathfrak{S}_{\sfd}}$. Moreover, there exists a representation $V^{\lambda}_w$ of $\emph{GL}_{\sfd}(\CC)$ such that $\emph{ch}(V^{\lambda}_w)=s^w_{\lambda}(x)$. 
\end{thm}
\begin{proof}
For each $\lambda\in\calP(\sfd)$, we can find in \cite[p.4433]{MolevSagan} the expression
\begin{equation}\label{factorial Schur in terms of Schur}
 s_{\lambda}(x|-a)=\sum_{\nu\subset\lambda}\bar g_{\lambda\nu}(a)s_{\nu}(x) 
\end{equation}
where
\begin{equation}\label{factorial Schur in terms of Schur2}
\bar g_{\lambda\nu}(a) =  \sum_{T\in\mathcal{T}(\lambda,\nu)}\prod_{\substack{\alpha\in\lambda \\ T(\alpha)\ \text{unbarred}}} a_{T(\alpha)+c(\alpha)}
\end{equation}
and $\mathcal{T}(\lambda,\nu)$ is a subset of semi-standard tableaux of shape $\lambda$ depending on $\nu$. For the completeness, we recall the definition of $\mathcal{T}(\lambda,\nu)$. Consider a sequence of partitions 
\begin{equation}\label{sequence of partitions}
 R : \emptyset=\rho^{(0)} \rightarrow \rho^{(1)} \rightarrow \cdots \rightarrow \rho^{(l)} = \nu.
\end{equation}
Let $r_i$ be the row number of the box that one added to $\rho^{(i-1)}$ to obtain $\rho^{(i)}$. An element of $\mathcal{T}(\lambda,R)$ is a semi-standard tableaux $T$ of shape $\lambda$ with entries in $\{1,\cdots,\sfd\}$, together with a sequence of boxes $\alpha_1,\cdots,\alpha_{l}$ in $\lambda$ such that the column order of $\alpha_i$ are strictly increasing and $T(\alpha_i)=r_i$. These entries $T(\alpha_i), i=1,\dots,l$ are called \emph{barred}. Let $\mathcal{T}(\lambda,\nu):= \bigsqcup_{R}\mathcal{T}(\lambda,R)$ where the union is taken over all sequence of the form (\ref{sequence of partitions}). 

By the definition (\ref{factorial Schur in terms of Schur2}), we see that $\bar g_{\lambda\nu}(a)$ is a homogeneous polynomial in $a$'s of degree $|\lambda|-|\nu|$ and with positive integral coefficients. By substituting $a_i \mapsto -(w_i/u)x_{\bar{\emptyset}}$ (assume $u=1$) to the equation (\ref{factorial Schur in terms of Schur}),  we obtain
\begin{equation}\label{char1}
s^{w}_{\lambda} (x)  = s_{\lambda}(x|(-w_1)x_{\bar{\emptyset}},(-w_2)x_{\bar{\emptyset}},\cdots)  = \sum_{\nu\subset\lambda} \bar g_{\lambda\nu}(w) s_{\divi}(x)^{|\lambda|-|\nu|}  s_{\nu}(x). 
\end{equation}
Define $K_{\nu, k}^{\mu}\in\ZZ_{\geq0}$ by the equation
\[
s_{\divi}(x)^k s_{\nu}(x) = \sum_{\mu \in \calP(\sfd)} K_{\nu,k}^{\mu}s_{\mu}(x),  \ \ \ \  k\in\ZZ_{\geq0}.
\]
Note that $K_{\nu,k}^{\mu}=0$ unless $\nu\subset\mu$ and $|\mu|=|\nu|+k$. From (\ref{char1}), we have
\[
s^{w}_{\lambda} (x) 
=\sum_{\mu \in \calP(\sfd)} \left( \sum_{\nu\subset\lambda}  \bar g_{\lambda\nu}(w) K_{\nu, |\lambda|-|\nu|}^{\mu}\right) s_{\mu}(x). 
\]
Observe that the coefficient of $s_{\mu}(x)$
\[
m_{\lambda\mu}(w):=\sum_{\nu\subset\lambda} \bar g_{\lambda,\nu}(w) K_{\nu, |\lambda|-|\nu|}^{\mu}
\]
is a non-negative integer and vanishes unless $|\mu|=|\lambda|$. Thus $s^{w}_{\lambda} (x)  \in \ZZ[x]^{\mathfrak{S}_{\sfd}}$ and it is the character of the representation
\[
V^{\lambda}_w:=\bigoplus_{\mu\in\calP(\sfd) \atop{|\lambda|=|\mu|}} (V^{\mu})^{\oplus m_{\lambda\mu}(w)}
\]
where $V^{\mu}$ is the irreducible representation of $\text{GL}_{\sfd}(\CC)$ with the highest weight $(\mu_1,\cdots,\mu_{\sfd})$.
\end{proof}
\begin{rem}
Since $s^w_{\lambda}(x)$ is a homogeneous polynomial, the representation $V^{\lambda}_w$ is homogeneous. That is, there exists a representation $S^{\lambda}_w$ of $\mathfrak{S}_{|\lambda|}$ such that 
\begin{align*}
 V^{\lambda}_w =
 (\CC^{\sfd})^{\otimes |\lambda|} \otimes_{\CC[\mathfrak{S}_{|\lambda|}]} S^{\lambda}_w
\end{align*}
as representations of $\text{GL}_{\sfd}(\CC)$ where $\GL_{\sfd}(\CC)$ naturally acts on the first factor in the right-hand-side.
\end{rem}
\begin{exm}
Suppose that $w_1=2, w_2=1, w_3=0$, and $u=1$. Then 
\begin{eqnarray*}
s^w_{\divi}(x)&=& 4 s_{\divi}(x), \\
 s^w_{(2)}(x) &=& 6 s_{(2)}(x) + 5 s_{(1,1)}(x).
\end{eqnarray*}
\end{exm}
Now we interpret the weighted Schubert structure constants in terms of representations. Suppose that $u=1$ and the weights are non-increasing ; $w_1\geq w_2\geq \cdots$. For all $\sfN>\sfd$, Corollary 5.6 in \cite{AbeMatsumura} shows that, the structure constants with weights $(w_{\sfN},\cdots,w_1)$ are non-negative rational numbers. Therefore, by Theorem \ref{twisted Schur goes weighted Schubert} and Theorem \ref{twisted Schur and rep}, we see that, for each $\lambda,\mu, \nu\in\calP(\sfd)$, there exist non-negative integers $m_{\lambda\mu}\in\ZZ_{\geq1}$ and $m_{\lambda\mu}^{\nu}\in\ZZ_{\geq0}$ such that 
\begin{align*}
 (V^{\lambda}_w \otimes V^{\mu}_w)^{\oplus m_{\lambda\mu}} = \bigoplus_{\nu\in\calP(\sfd)}  (V^{\nu}_w)^{\oplus m_{\lambda\mu}^{\nu}}
 \quad \text{ as representations of $\text{GL}_{\sfd}(\CC)$.}
\end{align*}
We can think of the fractions $m_{\lambda\mu}^{\nu}/m_{\lambda\mu}$ as \textit{rational multiplicity} of $V^{\nu}_w$'s in the tensor product of $V^{\lambda}_w \otimes V^{\mu}_w$ and they are nothing but the structure constants of the weighted Grassmannians.
\section{Determinantal Formulas}\label{Determinantal Formula}
In this section, we give two determinal formulas for the weighted Schubert classes $\wS_{\lambda} $. One is in terms of the special weighted Schubert classes and the other is in terms of the Chern classes of the tautological bundles. Recall that the weights $\{w_i\}$ has the reversed order, compared with the one in \cite{CortiReid} and \cite{AbeMatsumura}. 
\subsection{Known formulas}
Let $\calS \inc \calE \surj \calQ$ be the sequence of the tautological bundles of $\Gr(\sfd,\sfn)$ where $\calE=\Gr(\sfd,\sfn)\times\CC^{\sfn}$ is the trivial bundle. Let $\calF^{\ell}$ be the subbundle of $\calE$ defined by the coordinate plane generated by the \emph{last} $\ell$ coordinates. Then
\[
\sum_{r\geq 0}c^{\sfT}_r(\calF^{\ell} - \calS) =  \frac{\prod_{i=1}^{\ell} (1 - b_i)}{\prod_{j=1}^{\sfd}(1 - \bfx_j)}
\]
where $\bfx_1,\dots, \bfx_{\sfd}$ are the $\sfT$-equivariant Chern roots of the dual of $\calS$. The following formula is well-known (\textit{cf.} \cite[p.17] {Macdonald92}):
\[
\tS_{\lambda}=\det \big[c^{\sfT}_{\lambda_i-i+j}(\calF^{\lambda_i-i+\sfd} - \calS) \big]_{1\leq i\leq j\leq \sfd}.
\]
The equality $c^{\sfT}_r(\calF^{\ell} - \calS) = \sum_{k=0}^{r} h_{k}(b_{\ell+1}, \dots, b_{\sfn}) c^{\sfT}_{r-k}(\calQ)$ follows from
\[
\sum_{r\geq 0}c^{\sfT}_r(\calF^{\ell} - \calS)  = \frac{1}{\prod_{i=\ell +1}^{\sfn}(1 - b_i)} \sum_{r\geq 0}c^{\sfT}_r(\calE - \calS) = \frac{1}{\prod_{i=\ell +1}^{\sfn}(1 - b_i)} \sum_{r\geq 0}c^{\sfT}_r(\calQ).
\]
Therefore we also have
\begin{equation}\label{det tautological}
\tS_{\lambda}= \det \left[\sum_{k=0}^{\lambda_i-i+j} h_{{k}}(b_{\lambda_i-i+\sfd+1}, \dots, b_{\sfn}) c^{\sfT}_{{\lambda_i-i+j-k}}(\calQ) \right]_{1\leq i\leq j\leq \sfd}.
\end{equation}
On the other hand,  it is shown in \cite{LakshmibaiRaghavanSankaran} that
\begin{equation}\label{det special}
\tS_{\lambda}= \det \left[ \sum_{k=0}^{j-1}h_k(b_{\lambda_i-i+1+\sfd},\cdots,b_{\lambda_i-i+j-k+\sfd}) \tS_{(\lambda_i-i+j-k)} \right]_{1\leq i\leq j\leq \sfd}.
\end{equation}
When we specialize the above formulas to $b_k=0$ for all $k\in \NN$, we obtain the same determinant formula for the Schubert class in the ordinary cohomology of the Grassmannian.
\subsection{Determinantal formulas for the weighted Grassmannian}
Let $\calS_a \inc \calE_a \surj \calQ_a$ be the pullback of the sequence of tautological bundles along $\aPl^{\times}(\sfd,\sfn) \to \Gr(\sfd,\sfn)$. By quotienting the sequence by $\sfD_w$, we obtain the tautological sequence  $\calS_w \inc \calE_w \surj \calQ_w$ of orbifold vector bundles over $\wGr(\sfd,\sfn)$. Under the isomorphisms
\[
\xymatrix{
H_{\sfT}^*(\Gr(\sfd,\sfn)) \ar[r]_{\sff^*}^{\cong} & H_{\sfK}^*(\aPlt(\sfd,\sfn)) & H_{\sfT_{\!w}}^*(\wGr(\sfd,\sfn)) \ar[l]^{\sff_w^*}_{\cong}, 
}
\]
the equivariant Chern classes of those tautological bundles coincide. As we have seen, the equivariant Schubert classes also coincide. Therefore, the formulas (\ref{det tautological}, \ref{det special}) hold for both $H_{\sfK}(\aPlt(\sfd,\sfn))$ and $H_{\sfT_{\!w}}(\wGr(\sfd,\sfn))$. In particular, (\ref{det tautological}) gives the following in the non-equivariant setting.
\begin{thm}\label{weighted Giambelli}
For each $\lambda \in \calP(\sfd,\sfn)$, 
\[
\wS_{\lambda} = \det \left[\sum_{k=0}^{\lambda_i-i+j} \big(c_1(\calS_w)/u\big)^{k}h_{k}(w_{\lambda_i-i+\sfd+1}, \dots, w_{\sfn}) c_{\lambda_i-i+j-k}(\calQ_w) \right]_{1\leq i\leq j\leq \sfd}.
\]
\end{thm}
\begin{proof}
From (\ref{det tautological}), we have
\[
\wtS_{\lambda} = \det \left[\sum_{k=0}^{\lambda_i-i+j} \tilde{h}_{k,i}(b^w) c^{\sfT_{\!w}}_{\lambda_i-i+j-k}(\calQ_w) \right]_{1\leq i\leq j\leq \sfd}
\]
where 
\[
\tilde{h}_{k,i}(b^w)=h_{k}(b_{\lambda_i-i+\sfd+1}^w + (w_{\lambda_i-i+\sfd+1}/u)c_1^{\sfT_{\!w}}(\calS_w) , \dots, b_{\sfn}^w + (w_{\sfn}/u)c_1^{\sfT_{\!w}}(\calS_w)).
\]
Then we have $\tilde{h}_{k,i}(0) =  \big(c_1(\calS_w)/u\big)^{k}h_{k}(w_{\lambda_i-i+\sfd+1}, \dots, w_{\sfn})$.
\end{proof}
Similarly, by substituting $b_i \mapsto  -(w_i/w_{\bar{\emptyset}})\wS_{\divi}$ to (\ref{det special}),  one obtains the following.
\begin{thm}\label{wS in terms specials}
For each $\lambda\in\calP(\sfd,\sfn)$, we have
\[
 \wS_{\lambda}
 = \det\left[\sum_{k=0}^{j-1} h_k(w_{\lambda_i-i+1+\sfd},\cdots,w_{\lambda_i-i+j-k+\sfd}) \left(\frac{\wS_{\divi}}{-w_{{\bar{\emptyset}}}}\right)^{k}  \wS_{(\lambda_i-i+j-k)}\right]_{1\leq i,j\leq \sfd}
\]
\end{thm}
\begin{exm}
Let $(\sfd,\sfn)=(2,4)$ and $\lambda=(1,1)$. Assume $u=1$. Then we have
{\small \[
\wS_{(1,1)}
=\left[\begin{matrix}
c_{1}(\calQ_w) + c_1(\calS_w)h_{1}(w_3, w_4)&
c_{2}(\calQ_w)+ c_1(\calS_w)h_{1}(w_3, w_4) c_{1}(\calQ_w) + c_1(\calS_w)^{2}h_{2}(w_3, w_4)\\
1 &
c_{1}(\calQ_w) + c_1(\calS_w)h_{1}(w_2, w_3, w_4)
\end{matrix}\right].
\]}
\!\!\!
If $(w_1,w_2,w_3,w_4)=(0,2,1,0)$, then 
\[
\wS_{(1,1)}= \det 
\left[\begin{matrix}
c_1(\wcalQ) + c_1(\wcalS)	&	c_2(\wcalQ) + c_1(\wcalQ)c_1(\wcalS) +c_1(\wcalS)^2 \\
1 						&	c_1(\wcalQ) + 3c_1(\wcalS)
\end{matrix}\right].
\]
\end{exm}
\begin{exm}
Suppose $\lambda=(1,1) \in \calP(2,4)$.
\[
\wS_{\lambda}
=\det\left[
\begin{matrix}
\wS_{(1)}&
\wS_{(2)}+h_1(w_{3}) \left(-\wS_{\divi}/w_{{\bar{\emptyset}}}\right) \wS_{(1)}\\
1&
\wS_{(1)} +h_1(w_2) \left(-\wS_{\divi}/w_{{\bar{\emptyset}}}\right)
\end{matrix}
\right].
\]
For  $(w_1,w_2,w_3,w_4)=(0, 2, 1,0)$ and $u=1$, we have
\[
\wS_{(1,1)} = \det
\left[\begin{matrix}
\wS_{\divi} & \wS_{(2)} - \frac{1}{3}\wS_{\divi}^2 \\
1 & \frac{1}{3}\wS_{\divi} 
\end{matrix}\right]
\]
The reader can also check this equality directly by the Pieri rule (\ref{wPR}).
\end{exm}
\begin{rem}\label{rmk on difference of two Giambelli}
In the ordinary case, the above two propositions coincide since $S_{(r)} = c_r(\calQ)$. In general, they are different formulas, reflecting the fact that $\wS_{(r)} \not= c_r(\calQ_w)$.
\end{rem}
\subsection{A presentation of $H^*(\Gr(\sfd,\sfn))$}
We conclude by giving a presentation of the ordinary cohomology ring of the weighted Grassmannians $\wGr(\sfd,\sfn)$ over $\QQ$, in terms of Chern classes.  It is a well-known fact that $H_{\sfT}^*(\Gr(\sfd,\sfn))$ is generated by $c_i^{\sfT}(\calS)$ and $c_j^{\sfT}(\calQ)$ as a $\QQ[b]$-algebra and the relations are given by
\[
c^{\sfT}(\calS) c^{\sfT}(\calQ) = c^T(\calE), \ \ \ \textit{i.e.} \   \sum_{i=0}^rc^{\sfT}_i(\calS)c^{\sfT}_{r-i}(\calQ) = (-1)^re_r(b_1,\dots, b_n).
\]
Therefore, by 
\[
\sff^*(b_i)=\sff_w^*(b_i^w + (w_i/u)c_1^{\sfT_{\!w}}(\calS_w)),
\]
we have 
\[
\sum_{i=0}^rc^{\sfT_{\!w}}_i(\calS_w)c^{\sfT_{\!w}}_{r-i}(\calQ_w) = c_r^{\sfT_{\!w}}(\calE_w) = (-1)^re_r\left(b_1^w + (w_1/u)c_1^{\sfT_{\!w}}(\calS_w), \dots, b_{\sfn}^w + (w_{\sfn}/u)c_1^{\sfT_{\!w}}(\calS_w) \right).
\]
By setting $b^w_i=0$, we find that
\begin{equation}\label{relation-w}
\sum_{i=0}^rc_i(\calS_w)c_{r-i}(\calQ_w) = c_r(\calE_w) = (-1)^r(c_1(\calS_w)/u)^r e_r(w_1, \dots, w_{\sfn}).
\end{equation}
Thus $H^*(\wGr(\sfd,\sfn))$ is generated by $c_i(\calS_w)$ and $c_j(\calQ_w)$ with the relation (\ref{relation-w}). Therefore
\begin{prop} We have the following graded ring isomorphism
\[
H^*(\wGr(\sfd,\sfn)) \cong \frac{\QQ[c_1,\dots, c_{\sfd}, \bar c_1,\dots,\bar c_{{\sfn-\sfd}}]}{\big(\  \sum_{i=0}^rc_i \bar c_{r-i} = (-1)^r(c_1/u)^r e_r(w_1, \dots, w_{\sfn}), \ \ {r=1,\cdots,\sfn} \ \big)}
\]
where $c_0=\bar c_0 =1$, $c_i=0$ if $i>\sfd$ and $\bar c_j=0$ if $j > \sfn-\sfd$.
\end{prop}
\bibliography{references}{}

\begin{thebibliography}{10}

\bibitem{AbeMatsumura}
{\sc Abe, H., and Matsumura, T.}
\newblock {E}quivariant {C}ohomology of {W}eighted {G}rassmannians and
  {W}eighted {S}chubert {C}lasses.
\newblock {\em Int Math Res Notices}, doi: 10.1093/imrn/rnu003 (2014).

\bibitem{CortiReid}
{\sc Corti, A., and Reid, M.}
\newblock Weighted {G}rassmannians.
\newblock In {\em Algebraic geometry}. de Gruyter, Berlin, 2002, pp.~141--163.

\bibitem{FominGelfandPostnikov}
{\sc Fomin, S., Gelfand, S., and Postnikov, A.}
\newblock Quantum {S}chubert polynomials.
\newblock {\em J. Amer. Math. Soc. 10}, 3 (1997), 565--596.

\bibitem{Ikeda}
{\sc Ikeda, T.}
\newblock Schubert classes in the equivariant cohomology of the {L}agrangian
  {G}rassmannian.
\newblock {\em Adv. Math. 215}, 1 (2007), 1--23.

\bibitem{Ivanov97}
{\sc Ivanov, V.~N.}
\newblock The dimension of skew shifted {Y}oung diagrams, and projective
  characters of the infinite symmetric group.
\newblock {\em Zap. Nauchn. Sem. S.-Peterburg. Otdel. Mat. Inst. Steklov.
  (POMI) 240}, Teor. Predst. Din. Sist. Komb. i Algoritm. Metody. 2 (1997),
  115--135, 292--293.

\bibitem{Ivanov04}
{\sc Ivanov, V.~N.}
\newblock Interpolation analogues of {S}chur {$Q$}-functions.
\newblock {\em Zap. Nauchn. Sem. S.-Peterburg. Otdel. Mat. Inst. Steklov.
  (POMI) 307}, Teor. Predst. Din. Sist. Komb. i Algoritm. Metody. 10 (2004),
  99--119, 281--282.

\bibitem{KnutsonTao}
{\sc Knutson, A., and Tao, T.}
\newblock Puzzles and (equivariant) cohomology of {G}rassmannians.
\newblock {\em Duke Math. J. 119}, 2 (2003), 221--260.

\bibitem{LakshmibaiRaghavanSankaran}
{\sc Lakshmibai, V., Raghavan, K.~N., and Sankaran, P.}
\newblock Equivariant {G}iambelli and determinantal restriction formulas for
  the {G}rassmannian.
\newblock {\em Pure Appl. Math. Q. 2}, 3, Special Issue: In honor of Robert D.
  MacPherson. Part 1 (2006), 699--717.

\bibitem{LascouxSchutzenberger}
{\sc Lascoux, A., and Sch{\"u}tzenberger, M.-P.}
\newblock Polyn\^omes de {S}chubert.
\newblock {\em C. R. Acad. Sci. Paris S\'er. I Math. 294}, 13 (1982), 447--450.

\bibitem{Macdonald92}
{\sc Macdonald, I.~G.}
\newblock Schur functions: theme and variations.
\newblock In {\em S\'eminaire {L}otharingien de {C}ombinatoire
  ({S}aint-{N}abor, 1992)}, vol.~498 of {\em Publ. Inst. Rech. Math. Av.} Univ.
  Louis Pasteur, Strasbourg, 1992, pp.~5--39.

\bibitem{MolevSagan}
{\sc Molev, A.~I., and Sagan, B.~E.}
\newblock A {L}ittlewood-{R}ichardson rule for factorial {S}chur functions.
\newblock {\em Trans. Amer. Math. Soc. 351}, 11 (1999), 4429--4443.

\bibitem{Okounkov}
{\sc Okounkov, A.}
\newblock Quantum immanants and higher {C}apelli identities.
\newblock {\em Transform. Groups 1}, 1-2 (1996), 99--126.

\end{thebibliography}
\bibliographystyle{acm}
\end{document}